\documentclass[a4paper,10pt,reqno]{amsart}
\usepackage{latexsym, stmaryrd, mathtools}
\usepackage{ae}
\usepackage{tikz}
\usetikzlibrary{snakes}
\usetikzlibrary{matrix}

\DeclareMathAlphabet{\mathitbf}{OML}{cmm}{b}{it}

\newcommand{\NN}{{\mathbb N}} % Natural numbers
\newcommand{\emptyword}{\ensuremath{\epsilon}}
\newcommand{\sym}{\mathcal{S}} % the symmetric group
\newcommand{\B}{\mathcal{B}} % the beta trees
\newcommand{\barB}{\bar{\B}} % the indecomposable beta trees
\newcommand{\A}{\mathcal{A}} % the avoiders
\newcommand{\barA}{\bar{\A}} % the indecomposable avoiders
\newcommand{\slot}{\ensuremath{-}}

\theoremstyle{plain}
\newtheorem{theorem}{Theorem}
\newtheorem*{theorem*}{Theorem}
\newtheorem{corollary}[theorem]{Corollary}
\newtheorem*{corollary*}{Corollary}
\newtheorem{lemma}[theorem]{Lemma}
\newtheorem*{lemma*}{Lemma}

\newtheorem*{proposition*}{Proposition}
\newtheorem{conjecture}[theorem]{Conjecture}
\newtheorem*{conjecture*}{Conjecture}

\theoremstyle{definition}

\newtheorem*{definition*}{Definition}

\newtheorem*{example*}{Example}

\newtheorem*{problem*}{Problem}

\theoremstyle{remark}

\newtheorem*{remark*}{Remark}

\newcommand{\dd}{\text{-}}
\newcommand{\pc}{3\dd 1\dd 4\dd 2}
\newcommand{\pg}{2\dd 41\dd 3}

\newcommand{\lra}{\longrightarrow}

\DeclareMathOperator{\st}{\mathrm{std}} % the standardization op.
\DeclareMathOperator{\card}{\mathrm{card}}

\DeclareMathOperator{\des}{\mathrm{des}}

\DeclareMathOperator{\asc}{\mathrm{asc}}
\DeclareMathOperator{\lmin}{\mathrm{lmin}}
\DeclareMathOperator{\lmax}{\mathrm{lmax}}
\DeclareMathOperator{\rmin}{\mathrm{rmin}}
\DeclareMathOperator{\rmax}{\mathrm{rmax}}
\DeclareMathOperator{\comp}{\mathrm{comp}}

\DeclareMathOperator{\ldr}{\mathrm{ldr}}
\DeclareMathOperator{\lir}{\mathrm{lir}}

\DeclareMathOperator{\leaves}{\mathrm{leaves}}
\DeclareMathOperator{\internal}{\mathrm{int}}
\DeclareMathOperator{\sub}{\mathrm{sub}}
\renewcommand{\root}{\operatorname{root}}
\DeclareMathOperator{\lpath}{\mathrm{lpath}}
\DeclareMathOperator{\rpath}{\mathrm{rpath}}
\DeclareMathOperator{\lsub}{\mathrm{lsub}}
\DeclareMathOperator{\rsub}{\mathrm{rsub}}
\DeclareMathOperator{\stem}{\mathrm{stem}}
\DeclareMathOperator{\stemh}{\mathrm{gamma}}
\DeclareMathOperator{\stemhm}{\mathrm{beta}}

\newcommand{\btree}{\mbox{\ensuremath{\beta(1,0)}-tree}}
\newcommand{\btrees}{\btree s}

\newcommand{\cfill}{black!30}
\newcommand{\cfilll}{white}

\newcommand{\ns}{4pt}
\newcommand{\nodestyle}{\tikzstyle{every node} = [font=\footnotesize]}
\newcommand{\discstyle}{\tikzstyle{disc} = 
  [ circle,thin,fill=\cfilll,draw=black, minimum size=\ns, inner sep=0pt ] }
\newcommand{\style}{
  \nodestyle
  \discstyle
}

\newcommand{\scl}{0.45}

\newcommand{\leaf}{
  \begin{tikzpicture}[ scale=\scl, baseline=-2.5pt ]
    \style
    \node [disc] (r) at (0,0) {};
  \end{tikzpicture}
}
\newcommand{\edge}{  
  \begin{tikzpicture}[ xscale=0.4, yscale=0.3, baseline=1.9pt ]
    \style
    \node [disc] (r) at (0,1) {};
    \node [disc] (1) at (0,0) {};
    \draw (r)  -- (1);
  \end{tikzpicture}
}

\newcommand{\bbx}{  
  \begin{tikzpicture}[ scale=\scl, baseline=-2pt ]
    \style
    \node [disc] (r) at (0,1) {};
    \node [disc] (1) at (0,0) {};
    \draw (r) node[above=1pt] 
          {\ensuremath{1}} -- (1) node[below=1pt] {\ensuremath{1}};
  \end{tikzpicture}
}
\newcommand{\bbbx}{  
  \begin{tikzpicture}[ scale=\scl, baseline=(r11.base) ]
    \style
    \node [disc] (r)   at (0, 0) {};
    \node [disc] (r1)  at (0,-1) {};
    \node [disc] (r11) at (0,-2) {};
    \draw 
      (r)   node[left] {1} -- 
      (r1)  node[left] {1} -- 
      (r11) node[left] {1};
  \end{tikzpicture}
}
\newcommand{\bbbxx}{  
  \begin{tikzpicture}[ scale=\scl, baseline=(r2.base) ]
    \style
    \node [disc] (r)  at (   0, 0) {};
    \node [disc] (r1) at (-0.6,-1) {};
    \node [disc] (r2) at ( 0.6,-1) {};
    \draw 
      (r) node[above] {2} -- (r1) node[below] {1} -- 
      (r)                 -- (r2) node[below] {1};
  \end{tikzpicture}
}

\newcommand{\bxxx}{
  \style
  \node [disc] (r1)  at ( 0,-1) {};
  \node [disc] (r11) at (-1,-2) {};
  \node [disc] (r12) at ( 0,-2) {};
  \node [disc] (r13) at ( 1,-2) {};
  \draw 
  (r1) -- (r11) node[below=1pt] {1} 
  (r1) -- (r12) node[below=1pt] {1}
  (r1) -- (r13) node[below=1pt] {1};
}
\newcommand{\eeev}{
  \style
  \node [disc] (r)    at ( 0, 3 ) {};
  \node [disc] (r1)   at ( 0, 2 ) {};
  \node [disc] (r11)  at ( 0, 1 ) {};
  \node [disc] (r111) at (-0.8, 0.1 ) {};
  \node [disc] (r112) at ( 0, 0.1 ) {};
  \draw 
  (r11) -- (r111) node[below=1pt] {1}
  (r11) -- (r112) node[below=1pt] {1};
}

\newcommand{\exampleforest}{
    \path
    node [disc] (1)   at (-2, -1) {} 
    node [disc] (11)  at (-3, -2) {}
    node [disc] (12)  at (-1, -2) {}
    node [disc] (112) at (-2, -3) {}
    node [disc] (111) at (-4, -3) {}
    node [disc] (2)   at ( 2, -1) {} 
    node [disc] (21)  at ( 1, -2) {}
    node [disc] (22)  at ( 2, -2) {}
    node [disc] (23)  at ( 3, -2) {};

    \draw % edges and labels of the nodes
    (1)  node[above left=-1pt]  {1} -- (11)  node[above left=-1pt]  {2}
    (1)                      -- (12)  node[below=1pt] {1}
    (11)                     -- (111) node[below=1pt] {1}
    (11)                     -- (112) node[below=1pt] {1}
    (2) node[above right=-1pt] {3}  -- (21)  node[below=1pt] {1}
    (2)                      -- (22)  node[below=1pt] {1}
    (2)                      -- (23)  node[below=1pt] {1};
  }
\newcommand{\exampletree}{
  \node [disc] (r) at (0,0) {}; % root
  \exampleforest
  \draw
  (r) node[above=1pt] {4} -- (1)
  (r)                     -- (2);
}

\newcommand{\exampletreee}{
  \node [disc] (r) at (0,0)  {}; % root
  \node [disc] (x) at (0,-1) {}; % extra edge to the root
  \exampleforest
  \draw
  (r) node[above=1pt] {5} -- (1)
  (r)                     -- (x) node[below=1pt] {1}
  (r)                     -- (2);
}

\newcommand{\tri}{
  \filldraw[fill=\cfill, draw=\cdraw] 
  (0,0) -- (-2,-2) -- (2,-2) -- cycle;
  \path
  node [disc] (1)    at (0, 0) {} 
  node [disc] (11)   at (1,-1) {};
}

\newcommand{\extree}{
  \begin{center}
    \begin{tikzpicture}[ scale=\scl, baseline=0pt ]
      \discstyle
      \exampletree
    \end{tikzpicture}
  \end{center}
}

\title[$\beta(1,0)$-trees and nonseparable permutations]
      {Decompositions and statistics for \\ 
        $\beta(1,0)$-trees and nonseparable permutations}

\author[A. Claesson]{Anders Claesson}
\author[S. Kitaev]{Sergey Kitaev}
\author[E. Steingr\'{\i}msson]{Einar Steingr\'{\i}msson}
\thanks{The
  research presented here was supported by grant no. 060005012 from
  the Icelandic Research Fund}

\address{The Mathematics Institute, Reykjavik University, Kringlan 1,
  103 Reykjavik, Iceland} 

\begin{document}

\begin{abstract}
  The subject of pattern avoiding permutations has its roots in
  computer science, namely in the problem of sorting a permutation
  through a stack.  A formula for the number of permutations of length
  $n$ that can be sorted by passing it twice through a stack (where
  the letters on the stack have to be in increasing order) was
  conjectured by West, and later proved by Zeilberger.  Goulden and
  West found a bijection from such permutations to nonseparable planar
  maps, and later, Jacquard and Schaeffer presented a bijection
  from these planar maps to certain labeled plane trees, called
  \btrees. Using generating trees, Dulucq, Gire and West showed that
  nonseparable planar maps are equinumerous with permutations avoiding
  the (classical) pattern $2413$ and the barred pattern $41\bar{3}52$;
  they called these permutations \emph{nonseparable}. 

  We give a new bijection between \btrees\ and permutations avoiding the
  dashed patterns 3-1-4-2 and 2-41-3. These permutations can be seen
  to be exactly the reverse of nonseparable permutations.  Our
  bijection is built using decompositions of the permutations and the
  trees, and it translates seven statistics on the trees into
  statistics on the permutations.  Among the statistics involved are
  ascents, left-to-right minima and right-to-left maxima for the
  permutations, and leaves and the rightmost and leftmost paths for
  the trees.

  In connection with this we give a nontrivial involution on the
  \btrees, which specializes to an involution on unlabeled rooted
  plane trees, where it yields interesting results.

  Lastly, we conjecture the existence of a bijection between
  nonseparable permutations and two-stack sortable permutations
  preserving at least four permutation statistics.
\end{abstract}

\maketitle

\thispagestyle{empty}

\section{Introduction}

In Exercise~2.2.1.5 Knuth \cite{Kn69v1} asks the reader to ``Show that
it is possible to obtain a permutation $p_1p_2\dots p_n$ from $12\dots
n$ using a stack if and only if there are no indices $i<j<k$ such that
$p_j<p_k<p_i$.'' The equivalent inverse problem is to sort a
permutation (into increasing order) through a stack, and the
characterization in Knuth's exercise states that a permutation can be
sorted through a stack if and only if it \emph{avoids the pattern}
2-3-1. (A permutation avoids the pattern 2-3-1 if there are no indices
$i<j<k$ such that $p_k<p_i<p_j$.)

In his Ph.D. thesis West~\cite{We90} considered the problem of sorting
a permutation by passing it twice through a stack, where the letters
on the stack have to be in increasing order when read from top to
bottom. He conjectured that the number of permutations on
$[n]=\{1,\dots,n\}$ so sortable is
$a_n=2(3n)!/\big((2n+1)!(n+1)!\big)$. This was first proved by
Zeilberger~\cite{Ze92}, who found the functional equation
$$x^2F^3+x(2+3x)F^2+(1-14x+3x^2)F+x^2+11x-1=0
$$ for the generating function $F=\sum_n a_nx^n$ and then used
Lagrange's inversion formula to solve it. 

West also had noted that $a_n$ is the number of rooted nonseparable
planar maps on $n+1$ edges as enumerated by Brown and
Tutte~\cite{Br63, BrTu64, Tu63}. Based on this fact two bijective
proofs were later found. One by Dulucq et al.~\cite{DuGiGu98,
  DuGiWe96} who establish the correspondence using generating trees
and eight different families of permutations with forbidden
subsequences as intermediate sets. The other by Goulden and
West~\cite{GoWe96}: Using Raney paths they unearthed a recursive
structure on the permutations that parallels the recursive structure
that Brown and Tutte had found on the planar maps.

Jacquard and Schaeffer~\cite{JaGi98} presented a bijection from rooted
nonseparable planar maps to certain labeled plane trees. These trees
are the combinatorial structure that most transparently embodies the
recursive structure that Brown and Tutte had found on the planar maps,
and here is their definition: A \emph{\btree} is a rooted plane tree
labeled with positive integers such that
\begin{enumerate}
\item Leaves have label $1$.
\item The root has label equal to the sum of its children's labels.     
\item Any other node has label no greater than the sum of its
  children's labels.
\end{enumerate}
Below is an example of such a tree.\\[-4ex]
\extree

The main result (Theorem~\ref{thm_f}) of this paper is a one-to-one
correspondence between {\btrees} on $n$ edges and permutations of
length $n$ that avoid the two dashed (or generalized) patterns 3-1-4-2 and 2-41-3
(we call these permutations {\em avoiders}). That correspondence
``translates'' seven different statistics from \btrees\ to avoiders.  To
be more precise, we say that two vectors $(s_1,s_2,\ldots,s_k)$ and
$(t_1,t_2,\ldots,t_k)$ of statistics on sets $A$ and $B$,
respectively, {\em have the same distribution} if
$$
\sum_{a\in A}{x_1^{s_1(a)}x_2^{s_2(a)}\cdots x_k^{s_k(a)}} = 
\sum_{b\in B}{x_1^{t_1(b)}x_2^{t_2(b)}\cdots x_k^{t_k(b)}}.
$$ The bijection $f$ in Theorem \ref{thm_f} shows that the first
vector below has the same distribution on \btrees\ as the
second has on avoiders (see Section~\ref{prel} for definitions):
$$
\begin{array}{lllllllll}
  ( &\sub,   &\leaves, &\root, &\lpath, &\rpath, &\lsub, &\stemhm & ) \\
  ( & \comp, &1+\asc,  &\lmax, &\lmin,  &\rmax,  &\ldr,  &\lir    & )
\end{array}
$$

Bijections between two sets (of combinatorial objects) that take one
vector of statistics to another typically reveal that the two sets are
structurally similar, in a sense illuminated by the statistics in
question.  By doing exhaustive computations of large sets of
statistics on our avoiders, two-stack sortable permutations and
\btrees, respectively, we have found that there are many more, and
larger, vectors of statistics that are equidistributed between
avoiders and \btrees\ than between two-stack sortable permutations and
the trees.  Moreover, there are few and small equidistributions
between avoiders and two-stack sortables.  This suggests that the
avoiders are structurally more similar to the \btrees\ than two-stack
sortable permutations are.

Thus, it seems that the avoiders are a better set of permutations than
two-stack sortables to capture interesting properties of the trees
and, by extension, the nonseparable planar maps.  In fact, we were led
to some of the statistics on trees we study here---statistics that
belong to the vector of seven statistics mentioned above---by
well-known statistics on permutations. In particular, the left
subtrees presented themselves when we studied what the bijection $f$
translates the permutation statistic $\ldr$ to on \btrees. The left
subtrees are obtained by ``cutting'' the tree at each 1 on the left
path (except at the root). The statistic $\ldr$ is the place of the
first ascent in a permutation (if there is at least one ascent, and
equal to the length of the permutation otherwise).

In connection with Theorem~\ref{thm_f} mentioned above it should be
acknowledged that Dulucq et al.~\cite{DuGiGu98, DuGiWe96} gave a
bijection between rooted nonseparable planar maps and what they call
\emph{nonseparable permutations} (permutations avoiding the pattern
2-4-1-3 and the barred pattern 4-1-$\overline{3}$-5-2). It is not hard
to see that those permutations are the reverse of the permutations
avoiding 3-1-4-2 and 2-41-3. Dulucq et al. construct their bijection
via a generating tree of nonseparable permutations and a generating
tree of rooted nonseparable planar maps, and they show that their
bijection sends the pair $(\text{degree of root face},\,\text{number
  of nodes})$ to the pair $(\rmax, \des)$.

In a recent paper, Bonichon et al.\ ~\cite{BoBoEr08} give a
beautifully simple and direct bijection from Baxter permutations to
plane bipolar orientations.  The restriction of this bijection to
nonseparable permutations (a subset of the Baxter permutations) is the
same bijection to rooted nonseparable planar maps as the one given by
Dulucq et al.\ ~\cite{DuGiGu98}. Bonichon et al.\ also show that their
bijection translates five natural statistics on permutations to
statistics on the maps, thus strengthening the corresponding result of
Dulucq et al. A question then presents itself: If we take the
bijection given in the present paper and compose it with the bijection
of Jacquard and Schaeffer~\cite{JaGi98}, is the resulting bijection
equal to or different from the bijection of Dulucq et al?  The answer
is that it is different, and it is different also "up to symmetry."

Let us make this a bit more precise. It can be seen that of the seven
non-identity symmetries on permutations formed by composing reverse
($r$), complement ($c$) and inverse ($i$), the set of avoiders is
closed under three: $r\circ c$, $r\circ i$ and $c\circ i$.  We also
consider a symmetry on \btrees, namely mirror, which recursively
reverses the order of subtrees. On planar maps we consider two
symmetries: mirror and the planar dual (see \cite{BoBoEr08}). We claim
that our bijection composed with any of the mentioned symmetries gives
a bijection that is different from the bijection of Dulucq et
al.~\cite{DuGiGu98}. This can be seen by applying all these
bijections to the avoider $245316$.

Bousquet-M{\'e}lou~\cite{BoMe98} studied the generating function for
2-stack sortable permutations $\pi$ with respect to $5$ parameters:
length of $\pi$, $\des(\pi)$, $\lmax(\pi)$, $\rmax(\pi)$, and the
largest $i$ such that $(n, n-1, \dots, n+1-i)$ is a subsequence of
$\pi$ (where $n$ is the length of $\pi$). She showed that this
five-variable generating function is algebraic of degree 20.  As
mentioned above, we treat all these statistics, and more, on avoiders
and relate them to statistics on \btrees.

In Section~\ref{h} we introduce an involution $h$ on {\btrees}. It
gives us three further results about equidistributions (see
Theorem~\ref{thm_h}, Corollary~\ref{cor_hf},
Corollary~\ref{cor_hmf}). Moreover, $h$ gives an involution on
\emph{unlabeled} rooted plane trees.  Using that, we obtain some
interesting results on one-stack sortable permutations and also a
genuinely new bijection between (1-2-3)-avoiding and (1-3-2)-avoiding
permutations, yielding new equidistributions of statistics on these
two classes of permutations.  These results will be presented in a
forthcoming paper~\cite{CKS}.

\section{Preliminaries}\label{prel}

 Let $V=\{v_1,v_2,\ldots,v_n\}$ with $v_1<v_2<\dots<v_n$ be any finite
 subset of $\NN$. The \emph{standardization} of a permutation $\pi$ on
 $V$ is the permutation $\st(\pi)$ on $[n]$ obtained from $\pi$ by
 replacing the letter $v_i$ with the letter $i$. As an example,
 $\st(19452) = 15342$. If the set $V$ is fixed, the inverse of the
 standardization map is well defined, and we denote it by
 $\st^{-1}_V(\sigma)$; for instance, with $V=\{1,2,4,5,9\}$, we have
 $\st^{-1}_V(15342)=19452$. An \emph{occurrence} of the pattern {\pc}
 in a permutation $\pi=a_1a_2\dots a_n$ is a subsequence
 $o=a_ia_ja_ka_{\ell}$ (where $i<j<k<\ell$) of $\pi$ such that
 $\st(o)=3142$; an occurrence of {\pg} is a subsequence
 $o=a_ia_ja_{j+1}a_k$ such that $\st(o)=2413$. A permutation is said
 to \emph{avoid} a pattern if it has no occurrences of it. (For more
 on dashed permutation patterns see ~\cite{BaSt00,
 Cl01}.)\medskip

\emph{From now on we will refer to a permutation avoiding the two
  patterns {\pc} and {\pg} simply as an $\mathitbf{avoider}$}.\medskip

Here follows some more terminology that we shall use. An
\emph{interval} in a permutation is a factor (contiguous subsequence) that
contains a set of contiguous values. For example, $423$ is an interval
in $1642375$.  In particular, every permutation is an interval of
itself, and every letter is an interval too.

For words $\alpha$ and $\beta$ over the alphabet $\NN$ we define that
$\alpha\prec\beta$ if for all letters $a$ in $\alpha$ and all letters
$b$ in $\beta$ we have $a < b$. For instance, $412\prec 569$ and
$2\prec 348$. 
The reflexive closure of $\prec$ is a partial order on the set of
nonempty words over $\NN$. Indeed, if both $\alpha$ and $\beta$ are
nonempty then $\alpha\prec\beta$ if and only if $\max\alpha <
\min\beta$.

We now define the statistics on permutations we shall be concerned
with.  Let $\pi=a_1a_2\dots a_n$ be any permutation. An \emph{ascent}
is a letter followed by a larger letter; a \emph{descent} is a letter
followed by a smaller letter. The number of ascents and descents are
denoted $\asc(\pi)$ and $\des(\pi)$, respectively.  A
\emph{left-to-right minimum} of $\pi$ is a letter with no smaller
letter to the left of it; the number of left-to-right minima is
denoted $\lmin(\pi)$.  The statistics \emph{right-to-left minima}
($\rmin$), \emph{left-to-right maxima} ($\lmax$), and
\emph{right-to-left maxima} ($\rmax$) are defined similarly. The
statistic $\ldr(\pi)$ is defined as the largest integer $i$ such that
$a_1>a_2>\cdots>a_i$ (the leftmost decreasing run). Similarly,
$\lir(\pi)$ is defined as the largest integer $i$ such that
$a_1<a_2<\cdots<a_i$ (the leftmost increasing run). A \emph{component}
of $\pi$ is a nonempty factor $\tau$ of $\pi$ such that $\pi =
\sigma\tau\rho$ with $\sigma\prec\tau\prec\rho$, and such that if
$\tau$ = $\alpha\beta$ and $\alpha\prec\beta$ then $\alpha$ or $\beta$
is empty. By $\comp(\pi)$ we denote the number of components of
$\pi$. For instance, $\comp(213645)=3$, the components being $21$,
$3$, and $645$.

We also define some statistics on {\btrees}. By $\leaves(t)$ we denote
the number of leaves in $t$; by $\internal(t)$ we denote the number of
internal nodes (or nonleaves) in $t$.  Note that the root is an
internal node. By $\root(t)$ we denote the label of the root. The
number of subtrees (or, equivalently, the number of children of the
root) is denoted $\sub(t)$. The \emph{left-path} is the path from the
root to the leftmost leaf, and the \emph{right-path} is the path from
the root to the rightmost leaf.  The lengths of (number of edges on)
the left and right paths are denoted $\lpath(t)$ and $\rpath(t)$,
respectively. By $\stem(t)$ we denote the number of internal
nodes that are common to the left and the right-path. Also, $\lsub(t)$
and $\rsub(t)$ denote the number of $1$'s below the root on the left
and right paths, respectively.

We define $\stemhm(t)$ as follows: Order the leaves of a tree $t$ from
left to right and call them $\ell_1,\ell_2,\ldots,\ell_m$ (where
$\ell_1$ is leftmost and so on).  Look at the path from $\ell_1$ to
the root.  If no node on that path, except for the leaf $\ell_1$, has
label 1, reduce the labels on all nodes on that path by 1 and delete
$\ell_1$. Note that the resulting tree is a $\btree$ and that its
leaves are $\ell_2,\ldots,\ell_m$. Now look at $\ell_2$ and repeat
the process, until we come to a leaf $\ell_i$ whose path to the root
contains a node (other than $\ell_i$) that now has label 1. Then
$\stemhm(t)=i$. We end these preliminaries with an example: 

$$
t\; =\,
\begin{tikzpicture}[ scale=\scl, baseline=(r2).base ]
  \discstyle
  \path
  node [disc] (r)   at ( 0,  0.3) {} 
  node [disc] (r1)  at (-2, -1) {}
  node [disc] (r2)  at ( 0, -1) {}
  node [disc] (r21) at (-1, -2.1) {}
  node [disc] (r22) at ( 0, -2.1) {}
  node [disc] (r23) at ( 1, -2.1) {}
  node [disc] (r3)  at ( 2, -1) {}
  node [disc] (r31) at ( 2, -2.1) {};
  
  \draw % edges and labels of the nodes
  (r)  node[above=1pt] {4} -- (r1)  node[below=1pt]  {1}
  (r)                      -- (r2)  node[above left=-1pt] {2}
  (r2)                     -- (r21) node[below=1pt] {1}
  (r2)                     -- (r22) node[below=1pt] {1}
  (r2)                     -- (r23) node[below=1pt] {1}
  (r)                      -- (r3)  node[right=1pt] {1}
  (r3)                     -- (r31) node[below=1pt] {1};
\end{tikzpicture}
\qquad
\begin{array}{c}
  \ \\[-0.5ex]
  \leaves(t) = 5;\; \root(t) = 4; \\[1ex]
  \internal(t) = \sub(t) = \stemhm(t) = 3;\\[1ex]
  \rpath(t) = \rsub(t) = 2; \\[1ex]
  \lsub(t) = \lpath(t) = \stem(t) = 1.
\end{array}
$$

\section{The structure of {\btrees}}\label{btrees-structure}

We say a {\btree} on two or more nodes is \emph{indecomposable} if its
root has exactly one child and \emph{decomposable} if it has more than
one child. The {\btree} on one node is neither indecomposable nor
decomposable. Let $\B_n$ be the set of all {\btrees} on $n$ nodes,
and let $\barB_n$ be the subset of $\B_n$ consisting of the
indecomposable trees. Let $\B_n^k$ be the subset of $\B_n$ consisting
of the trees with root label $k$. For instance,
$$ \B_3=\big\{\;\,\bbbx\;,\;\bbbxx\;\big\}\qquad\;
\barB_3=\B_3^1=\big\{\;\bbbx\;\;\,\big\}\qquad\;
\B_3^2=\big\{\;\bbbxx\;\big\}
$$
Decomposable trees can be regarded as sums of indecomposable ones:
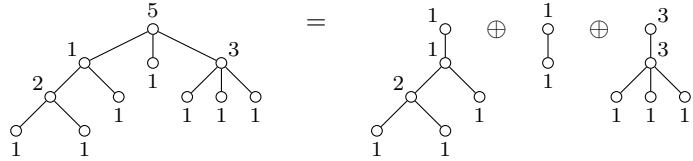
\begin{center}
  \begin{tikzpicture}[ scale=\scl, baseline=0pt ]
    \style    
    \exampletreee
  \end{tikzpicture}
  $\quad = \quad$
  \begin{tikzpicture}[ scale=\scl, baseline=0pt ]
    \style
    \path
    node [disc] (r1)  at (-2,  0) {} 
    node [disc] (rx)  at ( 1,  0) {} 
    node [disc] (r2)  at ( 4,  0) {}; 
%    \exampleforest
    \path
    node [disc] (1)   at (-2, -1) {} 
    node [disc] (11)  at (-3, -2) {}
    node [disc] (12)  at (-1, -2) {}
    node [disc] (112) at (-2, -3) {}
    node [disc] (111) at (-4, -3) {}
    node [disc] (x)   at ( 1, -1) {} 
    node [disc] (2)   at ( 4, -1) {} 
    node [disc] (21)  at ( 3, -2) {}
    node [disc] (22)  at ( 4, -2) {}
    node [disc] (23)  at ( 5, -2) {};

    \draw % edges and labels of the nodes
    (1)  node[above left=-1pt]  {1} -- (11)  node[above left=-1pt]  {2}
    (1)                      -- (12)  node[below=1pt] {1}
    (11)                     -- (111) node[below=1pt] {1}
    (11)                     -- (112) node[below=1pt] {1}
    (2) node[above right=-1pt] {3}  -- (21)  node[below=1pt] {1}
    (2)                      -- (22)  node[below=1pt] {1}
    (2)                      -- (23)  node[below=1pt] {1};

    \draw
    (r1) node[above left=-1pt]  {1} -- (1)   
    (rx) node[above=1pt]        {1} -- (x) node[below=1pt] {1}   
    (r2) node[above right=-1pt] {3} -- (2);
    \node[font=\normalsize] at (-0.5,0) {$\oplus$};
    \node[font=\normalsize] at (2.5,0) {$\oplus$};
  \end{tikzpicture}
\end{center}
\smallskip 
In fact we do not need to require $u$ and $v$ to be indecomposable for
the sum $u\oplus v$ to make sense. In general, we define that the root
label of $u\oplus v$ is the sum of the root label of $u$ and the root
label of $v$, and that the subtrees of $u\oplus v$ are those of $u$
followed by those of $v$. So,
$$
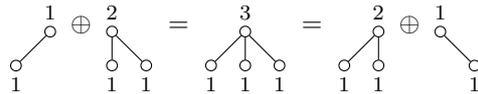
\begin{tikzpicture}[ scale=\scl, baseline=(r.base) ]
  \style
  \node [disc] (r) at (0,1) {};
  \node [disc] (1) at (-1,0) {};
  \draw (r) node[above=1pt] 
        {\ensuremath{1}} -- (1) node[below=1pt] {\ensuremath{1}};
\end{tikzpicture}\oplus
\begin{tikzpicture}[ scale=\scl, baseline=(r1.base) ]
  \style
  \node [disc] (r1)  at ( 0,-1) {};
  \node [disc] (r11) at ( 1,-2) {};
  \node [disc] (r12) at ( 0,-2) {};
  \draw 
  (r1) -- (r11) node[below=1pt] {1} 
  (r1) -- (r12) node[below=1pt] {1}; 
  \draw (r1) node[above=1pt] {2};
\end{tikzpicture}
=
\begin{tikzpicture}[ scale=\scl, baseline=(r1.base) ]
  \bxxx
  \draw (r1) node[above=1pt] {3};
\end{tikzpicture} 
=
\begin{tikzpicture}[ scale=\scl, baseline=(r1.base) ]
  \style
  \node [disc] (r1)  at ( 0,-1) {};
  \node [disc] (r11) at (-1,-2) {};
  \node [disc] (r12) at ( 0,-2) {};
  \draw 
  (r1) -- (r11) node[below=1pt] {1} 
  (r1) -- (r12) node[below=1pt] {1}; 
  \draw (r1) node[above=1pt] {2};
\end{tikzpicture} \oplus
\begin{tikzpicture}[ scale=\scl, baseline=(r.base) ]
  \style
  \node [disc] (r) at (0,1) {};
  \node [disc] (1) at (1,0) {};
  \draw (r) node[above=1pt] 
        {\ensuremath{1}} -- (1) node[below=1pt] {\ensuremath{1}};
\end{tikzpicture}
$$ 

Further, there is a simple one-to-one correspondence $\lambda$ between
the Cartesian product $[k]\times\B_{n-1}^k$ and the disjoint union
$\cup_{i=1}^k\barB_{n}^i$, where $\barB_n^k$ is the subset of
$\barB_n$ consisting of the trees with root label $k$:
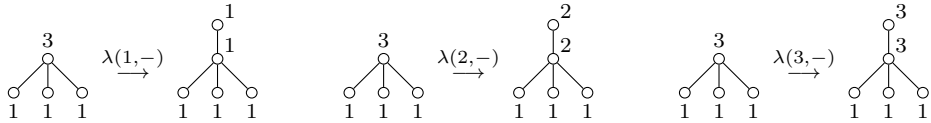
$$
\begin{tikzpicture}[ scale=\scl, baseline=(r11.base) ]
  \bxxx
  \draw (r1) node[above=1pt] {3};
\end{tikzpicture}
\,\raisebox{2ex}{$\substack{\lambda(1,\slot)\\ \lra}$}\,
\begin{tikzpicture}[ scale=\scl, baseline=(r11.base) ]
  \bxxx
  \node [disc] (r)  at (0,0) {};
  \draw (r) node[above right=-1pt] {1} -- (r1) node[above right=-1pt] {1};
\end{tikzpicture}
\qquad\;\;
\begin{tikzpicture}[ scale=\scl, baseline=(r11.base) ]
  \bxxx
  \draw (r1) node[above=1pt] {3};
\end{tikzpicture}
\,\raisebox{2ex}{$\substack{\lambda(2,\slot)\\ \lra}$}\,
\begin{tikzpicture}[ scale=\scl, baseline=(r11.base) ]
  \bxxx
  \node [disc] (r)  at (0,0) {};
  \draw (r) node[above right=-1pt] {2} -- (r1) node[above right=-1pt] {2};
\end{tikzpicture}
\qquad\;\;
\begin{tikzpicture}[ scale=\scl, baseline=(r11.base) ]
  \bxxx
  \draw (r1) node[above=1pt] {3};
\end{tikzpicture}
\,\raisebox{2ex}{$\substack{\lambda(3,\slot)\\ \lra}$}\,
\begin{tikzpicture}[ scale=\scl, baseline=(r11.base) ]
  \bxxx
  \node [disc] (r)  at (0,0) {};
  \draw (r) node[above right=-1pt] {3} -- (r1) node[above right=-1pt] {3};
\end{tikzpicture}
$$ In general, if $t$ is a tree with root label $k$ and $i$ is an
integer such that $1\leq i\leq k$, then $\lambda(i,t)$ is obtained
from $t$ by joining a new root via an edge to the old root; and both
the new root and the old root are assigned the label $i$.

Thus each \btree, $t$, is of exactly one the following three forms:
\begin{itemize}
\item[] $t=\leaf$,      \hfill (the one node tree)
\item[] $t=u\oplus v$,   \hfill (decomposable)
\item[] $t=\lambda(i,u)$, where $1\leq i\leq \root u$,\hfill(indecomposable)
\end{itemize}
in which $u$ and $v$ are \btrees. Note that any tree that is
decomposable with respect to $\oplus$ (second case above) is
indecomposable with respect to $\lambda$; that is, it is not of the
form $\lambda(i,u)$. Also, any tree that is indecomposable with
respect to $\oplus$ (third case above) is decomposable with respect to
$\lambda$. Thus each tree with at least one edge is decomposable with
respect to exactly one of $\oplus$ or $\lambda$, and if we keep
decomposing until only single node trees remain we get an unambiguous
encoding of \btrees, an example of which is given by
$$    
\begin{tikzpicture}[ scale=\scl, baseline=(111.base) ]
  \style 
  \path
  node [disc] (1)   at ( 0, -1) {} 
  node [disc] (11)  at ( 0, -2) {}
  node [disc] (111) at (-.65, -3) {}
  node [disc] (112) at (0.65, -3) {};
  
  \draw % edges and labels of the nodes
  (1) node[right=2pt] {2} -- (11) node[right=2pt] {2}
  (11)                   -- (111) node[below=1pt]     {1}
  (11)                   -- (112) node[below=1pt]     {1};
\end{tikzpicture}
= \lambda\Big(2,
\begin{tikzpicture}[ scale=\scl, baseline=(111.base) ]
  \style 
  \path
  node [disc] (11)  at ( 0, -2) {}
  node [disc] (111) at (-0.65, -3) {}
  node [disc] (112) at ( 0.65, -3) {};
  
  \draw % edges and labels of the nodes
  (11) node[above=2pt] {2} -- (111) node[below=1pt] {1}
  (11)                     -- (112) node[below=1pt] {1};
\end{tikzpicture}
\Big)
= \lambda\Big(2,\,\bbx\oplus\bbx\,\Big)
= \lambda\Big(2,\,\lambda(1,\leaf)\oplus\lambda(1,\leaf)\,\Big).
$$

\section{The structure of avoiders}

In this paper we construct a bijection between avoiders and {\btrees}
by defining a sum on permutations analogous to the sum on trees, and a
function $\phi$ analogous to $\lambda$; the sum on permutations works
like this:
$$21 \oplus 132 = 21354.
$$ In general, $\sigma\oplus\tau = \sigma\tau'$ where $\tau'$ is
obtained from $\tau$ by adding $|\sigma|$ to each of its letters. We
call a nonempty permutation $\pi$ \emph{decomposable} if it can be
written $\pi=\sigma\oplus\tau$ with $\sigma$ and $\tau$ both nonempty;
otherwise we call it \emph{indecomposable}. 

Describing how the map $\phi$ works is quite a bit
harder. In the case of {\btrees}, $\lambda(i,t)$ is indecomposable and
has root label $i$. In the case of avoiders, $\phi(i,\pi)$ is
indecomposable and has $i$ left-to-right maxima. 
For instance,
$$\pi = \underline{2}1\underline{58}6473
$$ has $3$ left-to-right maxima, namely the underlined elements, $2$,
$5$, and $8$. The $3$ images of $\pi$ under the function $\phi$ are
$\phi(1,\pi)$, $\phi(2,\pi)$, and $\phi(3,\pi)$; they should have $1$,
$2$, and $3$ left-to-right maxima, respectively. To achieve this we
start by inserting $9$ immediately before the $i$th left-to-right
maximum:
$$
\pi_1=\underline{9}21586473;\; \quad
\pi_2=\underline{2}1\underline{9}586473;\; \quad
\pi_3=\underline{2}1\underline{59}86473.
$$ Clearly $\pi_i$ has $i$ left-to-right maxima. Note also that it is
necessary to insert $9$ \emph{immediately} before the $i$th
left-to-right maximum, or else a $\pg$ pattern would be formed. To be
concrete, there is only one position in $\pi$ between the first
left-to-right maximum, $2$, and the second left-to-right maximum, $5$,
where we can insert $9$. Placing $9$ before $1$ would yield
$291586473$ and, as witnessed by the subsequence $2915$, that
permutation contains the pattern $\pg$. In general, inserting $n$
between two consecutive left-to-right maxima, say $b$ and $c$ with
$b<c$, but not immediately before $c$, will lead to an occurrence of
\pg\ formed by $b$, $n$, the element immediately to the right of $n$,
and $c$.

A permutation obtained by inserting a new largest element in an
avoider as described above will have the correct number of
left-to-right maxima; it will also avoid $\pc$ and $\pg$; it will
however \emph{not} be indecomposable, in general. Above, $\pi_1$ is
indecomposable but $\pi_2$ and $\pi_3$ are not. There is in fact a
simple criterion for indecomposability of avoiders: If $\pi$ is a
permutation in which the largest letter precedes the smallest letter,
then clearly $\pi$ is indecomposable. For $(\pc)$-avoiding
permutations the converse is also true.

\begin{lemma}\label{n-precedes-1}
  In any indecomposable $(\pc)$-avoiding permutation the largest
  letter precedes the smallest letter, when read from left to right.
\end{lemma} 

\begin{proof}
  We shall demonstrate the contrapositive statement: if
  $\pi\in\sym_n$ and $1$ precedes $n$ in $\pi$, then either $\pi$
  contains an occurrence of $\pc$ or $\pi$ is decomposable. To this
  end, let a permutation $\pi=\sigma n\tau$ in $\sym_n$ with
  $1\in\sigma$ be given. If $\tau$ is empty, then $n$ is a component
  and thus $\pi$ is decomposable. If $\sigma \prec \tau$ then $n\tau$
  is a component; otherwise the subword
  $\rho=(\,x\in\sigma:x>\min\tau\,)$ of $\sigma$ is nonempty, and
  there are two possibilities: either $\rho$ is a factor in $\sigma$
  and $\pi=\sigma'\rho n\tau$ with $1\in\sigma'$; or there is an $x$
  in $\rho$ to the left of a $y$ in $\sigma$ such that
  $y<\min\tau$. In the former case, $\rho n\tau$ is a component; in
  the latter case, $(x,y,n,\min\tau)$ is an occurrence of the pattern
  $\pc$.
\end{proof}

We shall now describe a function $\psi$ that can turn a decomposable
avoider into an indecomposable one, while preserving many of its other
properties. For now, let us concentrate on $\pi_3=215986473$ (from
above), and let us describe how to get $\psi(\pi_3)$. In doing so it
will be convenient to refer to the following picture.

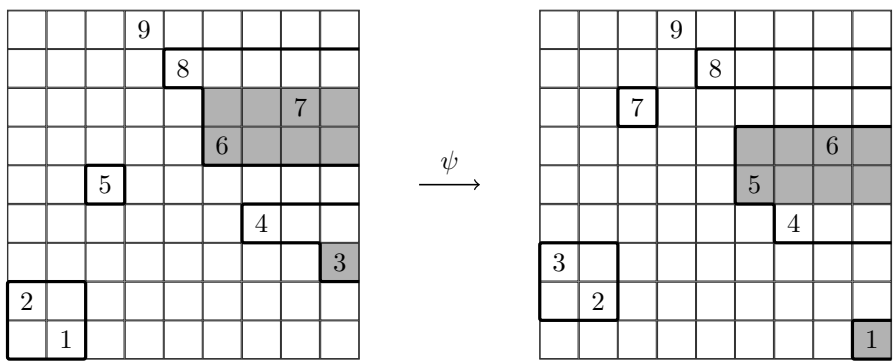
\begin{center}
  \begin{tikzpicture}
    \matrix at (0,0) [nodes={thin, draw=black!80, minimum size=5mm }] {
      \node (19) {};  & \node (29) {};  & \node (39) {};  & \node (49) {9}; & \node (59) {};  & \node (69) {};  & \node (79) {};  & \node (89) {};  & \node (99) {};  \\
      \node (18) {};  & \node (28) {};  & \node (38) {};  & \node (48) {};  & \node (58) {8}; & \node (68) {};  & \node (78) {};  & \node (88) {};  & \node (98) {};  \\
      \node (17) {};  & \node (27) {};  & \node (37) {};  & \node (47) {};  & \node (57) {};  & \node[fill=\cfill] (67) {};  & \node[fill=\cfill] (77) {};  & \node[fill=\cfill] (87) {7}; & \node[fill=\cfill] (97) {};  \\
      \node (16) {};  & \node (26) {};  & \node (36) {};  & \node (46) {};  & \node (56) {};  & \node[fill=\cfill] (66) {6}; & \node[fill=\cfill] (76) {};  & \node[fill=\cfill] (86) {};  & \node[fill=\cfill] (96) {};  \\
      \node (15) {};  & \node (25) {};  & \node (35) {5}; & \node (45) {};  & \node (55) {};  & \node (65) {};  & \node (75) {};  & \node (85) {};  & \node (95) {};  \\
      \node (14) {};  & \node (24) {};  & \node (34) {};  & \node (44) {};  & \node (54) {};  & \node (64) {};  & \node (74) {4}; & \node (84) {};  & \node (94) {};  \\
      \node (13) {};  & \node (23) {};  & \node (33) {};  & \node (43) {};  & \node (53) {};  & \node (63) {};  & \node (73) {};  & \node (83) {};  & \node[fill=\cfill] (93) {3}; \\
      \node (12) {2}; & \node (22) {};  & \node (32) {};  & \node (42) {};  & \node (52) {};  & \node (62) {};  & \node (72) {};  & \node (82) {};  & \node (92) {};  \\
      \node (11) {};  & \node (21) {1}; & \node (31) {};  & \node (41) {};  & \node (51) {};  & \node (61) {};  & \node (71) {};  & \node (81) {};  & \node (91) {};  \\
    };
    \draw [black!80] (19.north west) rectangle (91.south east);
    \draw [very thick, rounded corners=0.8pt] (12.north west) rectangle (21.south east);
    \draw [very thick, rounded corners=0.8pt] (35.north west) rectangle (35.south east);
    \draw [very thick, rounded corners=0.8pt] (98.north east) -- (58.north west) -- (58.south west) -- (58.south east) -- (66.south west)-- (96.south east);
    \draw [very thick, rounded corners=0.8pt] (94.north east) -- (74.north west) -- (74.south west) -- (93.north west) -- (93.south west) -- (93.south east);
    
    \draw[->, semithick] (3.1,0) -- (3.9,0) node[midway, yshift=3mm] {$\psi$};
    
    \matrix at (7,0) [nodes={thin, draw=black!80, minimum size=5mm }] {
      \node (19) {};  & \node (29) {};  & \node (39) {};  & \node (49) {9}; & \node (59) {};  & \node (69) {};  & \node (79) {};  & \node (89) {};  & \node (99) {};  \\
      \node (18) {};  & \node (28) {};  & \node (38) {};  & \node (48) {};  & \node (58) {8}; & \node (68) {};  & \node (78) {};  & \node (88) {};  & \node (98) {};  \\
      \node (17) {};  & \node (27) {};  & \node (37) {7}; & \node (47) {};  & \node (57) {};  & \node (67) {};  & \node (77) {};  & \node (87) {};  & \node (97) {};  \\
      \node (16) {};  & \node (26) {};  & \node (36) {};  & \node (46) {};  & \node (56) {};  & \node[fill=\cfill] (66) {};  & \node[fill=\cfill] (76) {};  & \node[fill=\cfill] (86) {6}; & \node[fill=\cfill] (96) {};  \\
      \node (15) {};  & \node (25) {};  & \node (35) {};  & \node (45) {};  & \node (55) {};  & \node[fill=\cfill] (65) {5}; & \node[fill=\cfill] (75) {};  & \node[fill=\cfill] (85) {};  & \node[fill=\cfill] (95) {};  \\
      \node (14) {};  & \node (24) {};  & \node (34) {};  & \node (44) {};  & \node (54) {};  & \node (64) {};  & \node (74) {4}; & \node (84) {};  & \node (94) {};  \\
      \node (13) {3}; & \node (23) {};  & \node (33) {};  & \node (43) {};  & \node (53) {};  & \node (63) {};  & \node (73) {};  & \node (83) {};  & \node (93) {}; \\
      \node (12) {};  & \node (22) {2}; & \node (32) {};  & \node (42) {};  & \node (52) {};  & \node (62) {};  & \node (72) {};  & \node (82) {};  & \node (92) {};  \\
      \node (11) {};  & \node (21) {};  & \node (31) {};  & \node (41) {};  & \node (51) {};  & \node (61) {};  & \node (71) {};  & \node (81) {};  & \node[fill=\cfill] (91) {1};  \\
    };
    \draw [black!80] (19.north west) rectangle (91.south east);
    \draw [very thick, rounded corners=0.8pt] (13.north west) rectangle (22.south east);
    \draw [very thick, rounded corners=0.8pt] (37.north west) rectangle (37.south east);
    \draw [very thick, rounded corners=0.8pt] (98.north east) -- (58.north west) -- (58.south west) -- (98.south east);
    \draw [very thick, rounded corners=0.8pt] (96.north east) -- (66.north west) -- (65.south west) -- (65.south east) -- (74.south west) -- (94.south east);
    \draw [very thick, rounded corners=0.8pt] (91.north east) -- (91.north west) -- (91.south west) -- (91.south east);
  \end{tikzpicture}  
\end{center}

To the left of $9$, the largest element in $\pi_3$, we have
$\sigma=215$; to the right we have $\tau=86473$. The patterns
(standardizations) of the left and right parts are $\bar\sigma=213$
and $\bar\tau=53241$. We shall make a new permutation
$\pi_3'=\psi(\pi_3)$ that will be similar to $\pi_3$ in the sense that
the letters to the left and right of $9$ in $\pi_3'$ also will form
the patterns $\bar\sigma$ and $\bar\tau$. Thus, to specify $\pi_3'$ it
is sufficient to give the $3$ element set, call it $L$, from which the
left part of $\pi_3'$ is built. Further, due to
Lemma~\ref{n-precedes-1}, $\pi_3'$ will be indecomposable if $1$ is
not a member of $L$. The underlying set of $\sigma$ is
$$\{1,2,5\}.
$$ 
Considering $\sigma=215$ as a permutation of this set, we can divide it
into two intervals $21$ and $5$. To each of the letters in the first
interval we will add some positive number $m_1+1$. Similarly, to the
letters in the second interval we will add some positive number
$m_2+1$. (In the picture above, $m_1+1$ and $m_2+1$ are the number of
elements in the gray areas.) The resulting underlying set of the left
part of $\pi_3'$ is
$$L=\big\{\,1+(m_1+1),\,2+(m_1+1),\,5+(m_2+1)\,\big\}.
$$ That $m_1+1$ is positive implies that $1$ is not in $L$,
which in turn implies that $\pi_3'$ is indecomposable. 

We now describe how $m_1$ and $m_2$ are determined. Let $w_1$ be
the subsequence of $\tau$ whose elements are between $21$ and $5$ in
value. So, $w_1=43$. Similarly, let $w_2$ be the subsequence of $\tau$
whose elements are larger than $5$ in value. So, $w_2=867$. 

Denote by $\hat 0_i$ the smallest element of $w_i$. Then $m_i$ is
the number of elements to the right of $\hat 0_i$ in $w_i$ that are
smaller than all elements to the left of $\hat 0_i$ in $w_i$.

We find that $m_1=0$ and $m_2=1$. Consequently, $L=\{2,3,7\}$. The
final step is to fill in the left part using the elements of $L$ and
to fill in the right part using the remaining elements $[8]\setminus
L=\{1,4,5,6,8\}$, while preserving the patterns $\bar\sigma$ and
$\bar\tau$. The resulting permutation is $327985461$.

We shall now go through the definition of $\psi$ again, this time
dealing with the general case. Let $\pi$ be an avoider whose first
letter is not $n$. In addition, assume that $n$ precedes $n-1$ in
$\pi$. Let $\sigma$ and $\tau$ be defined by $\pi=\sigma n\tau$. Note
that, by assumption, both $\sigma$ and $\tau$ are then nonempty. Let
$\sigma_1$, \dots, $\sigma_k$ be the intervals of $\sigma$, and, by
convention, let $\sigma_{0}=\emptyword$ and
$\sigma_{k+1}=\emptyword$. For each $i=0,\dots,k$ define a subsequence
$w_i$ of $\tau$ by
$$ w_i = (\,x\in\tau\,:\,\sigma_i\prec x \prec \sigma_{i+1}\,).
$$ Denote by $u_i$ and $v_i$ the parts of $w_i$ that are to the left
and to the right of the smallest element of $w_i$. In short, $w_i =
u_i\hat{0} v_i$ where $\hat 0_i=\min w_i$. We shall now specify the
set $L$ of elements to the left of $n$ in $\psi(\pi)$; this is the
key object in the definition of $\psi$: let
$m_i = \card\{\, x\in v_i: x \prec u_i\,\}$ and 
$L_i = \{\, x+m_i+1 :  x\in\sigma_i\,\}$; then $L = \cup_{i=1}^k L_i$.

Finally we define $\psi(\pi)$ as the result of filling in the
elements of $L$ to the left of $n$ while respecting the pattern
$\bar\sigma=\st(\sigma)$, and filling in the elements of
$R=[n-1]\setminus L$ to the right of $n$ while respecting the pattern
$\bar\tau=\st(\tau)$:
\begin{equation}\label{def_psi}
  \psi(\pi) = \sigma' n\tau',\text{ where }
  \sigma' = \st_L^{-1}(\bar\sigma) \text{ and }
  \tau' = \st_R^{-1}(\bar\tau).
\end{equation}
For the definition of $\st$, see Section \ref{prel}.

\begin{lemma}\label{lemma_psi0}
  Let $n\geq 2$, and let $\tilde\A_n$ be the set of avoiders whose
  first letter is not $n$. Then the function $\psi$, as defined in the preceding
  paragraph, is a bijection from 
  $$
  \big\{\,\pi\in\tilde\A_n\mid\text{$n$ precedes $n-1$ in $\pi$}\,\big\}
  \;\text{ onto }\; 
  \big\{\,\pi\in\tilde\A_n\mid\text{$\pi$ is indecomposable}\,\big\}.
  $$ 
\end{lemma}

\begin{proof}
  We shall use the same notation as above, so $\pi=\sigma n \tau$ is
  an avoider with $\sigma$ nonempty and the letter $n-1$ belongs to
  $\tau$. Further, we shall split this proof into 4 parts, showing
  that:
  \begin{itemize}
    \item[{\it a.}] the permutation $\psi(\pi)$ is indecomposable;
    \item[{\it b.}] the letter $n$ is not the first letter in $\psi(\pi)$;
    \item[{\it c.}] the permutation $\psi(\pi)$ avoids \pc\ and \pg;
    \item[{\it d.}] the function $\psi$ is injective;
    \item[{\it e.}] the function $\psi$ is surjective.
  \end{itemize}\smallskip

  {\it Part a.} Looking at the definition of $\psi$ we see that all the
  $m_i$s are nonnegative; hence all members of $L$ (i.e., all elements
  to the left of $n$ in $\psi(\pi)$) are bigger than $1$ and it
  follows by Lemma~\ref{n-precedes-1} that $\psi(\pi)$ is
  indecomposable.\medskip
 
  {\it Part b.} By assumption $\sigma$ is nonempty; by definition
  $\psi$ preserves the position of $n$. Thus $n$ is not the first
  letter in $\psi(\pi)$.\medskip

  {\it Part c.} We shall show that $\psi(\pi)=\sigma' n\tau'$ avoids
  $\pc$ and $\pg$, given that $\pi=\sigma n\tau$ does. Consider the
  contrapositive statement, and assume that there is an occurrence $o$
  of $\pc$ or $\pg$ in $\psi(\pi)$. Since the $L_i$s are intervals and
  the underlying permutations of the patterns $\pc$ and $\pg$ do not
  contain any nontrivial intervals we know that either $o$ is entirely
  contained in $\sigma'$ or only its first letter belongs to
  $\sigma'$. For the first case it suffices to recall that $\psi$
  preserves the patterns of $\sigma$ and $\tau$. The second case is
  more intricate:
  Let $o=c'a'd'b'$ be any occurrence of $\pc$ in $\psi(\pi)$ in which
  $c'$ is a letter of $\sigma'$ and $a'd'b'$ is a subword of $\tau'$.
  Let $a$, $b$, $c$, and $d$ be the preimages of $a'$, $b'$, $c'$, and
  $d'$ under $\psi$. Since $\psi$ preserves the pattern of $\tau$ we
  know that $a<b<d$. Also, $\tau'$ is the disjoint union of the
  sequences $\psi(w_i)$; let us call these the \emph{blocks} of
  $\tau'$. The assumption that $o$ is an occurrence of $\pc$ implies
  that $d'$ belongs to a different block than $a'$ and $b'$. This, in
  turn, implies that there is a letter $x$ in $\tau$ such that $b<x<d$
  and $x$ is to the left of $a$. Thus $x a d b$ is an occurrence of
  $\pc$ in $\pi$.

  Let $o=b'd'a'c'$ be any occurrence of $\pg$ in $\psi(\pi)$ in which
  $b'$ is a letter of $\sigma'$ and $d'a'c'$ is a subword of $\tau'$.
  As before, let $a$, $b$, $c$, and $d$ be the preimages of $a'$,
  $b'$, $c'$, and $d'$ under $\psi$. Since $\psi$ preserves the
  pattern of $\tau$ we have $a<c<d$. The assumption that $o$ is an
  occurrence of $\pg$ implies that $a'$ belongs to a different block
  than $c'$ and $d'$. If $a'=a$ then $b d a c$ is an occurrence of
  $\pg$ in $\pi$. If $a' \neq a$ then there is a letter $x$
  in $\tau$ such that $a<x<c$ and $x$ is to the left of $d$, and so $xd
  a c$ is an occurrence of $\pg$ in $\pi$.

  {\it Part d.} Note that the smallest letter in $w_1$ is to the left
  of any letter in $w_0$; otherwise, an occurrence $bdac$ of $\pg$ is
  materialized by letting $a=\psi(\pi)(i)$ be the first letter of
  $w_0$; $b$ be any letter in $\sigma_1$; $c$ be the smallest letter
  in $w_1$; and $d=\psi(\pi)(i-1)$ be the left neighbor of $a$. The
  following picture illustrates this argument in the special case when
  $\sigma$ is an interval.
  $$
  
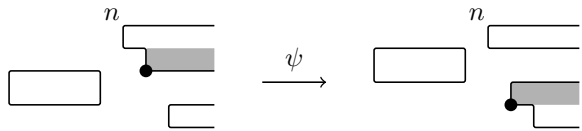
\begin{tikzpicture}[ scale = 0.3 ]
    \node at (4.5,2.5) {$n$};
    \draw[semithick, rounded corners=0.8pt] (0,-1.5) rectangle (4,0);
    \fill[fill=\cfill, rounded corners=0.8pt] (9,0) -- (6,0) -- (6,1) -- (9,1);
    \draw[semithick, rounded corners=0.8pt] (9,0) -- (6,0) -- (6,1) -- (5,1) -- (5,2) -- (9,2);
    \fill[fill=black] (6,0) circle (8pt);
    \draw[semithick, rounded corners=0.8pt] 
    (9,-2.5) -- (7,-2.5) -- (7,-1.5) -- (9,-1.5);
    \draw[->, semithick] (11.1,-.5) -- (13.9,-.5) 
    node[midway, yshift=3mm] {$\psi$};
    \node at (20.5,2.5) {$n$};
    \fill[fill=\cfill, rounded corners=0.8pt] 
    (25,-1.5) -- (22,-1.5) -- (22,-.5) -- (25,-.5);
    \draw[semithick, rounded corners=0.8pt] (16,-.5) rectangle (20,1);
    \draw[semithick, rounded corners=0.8pt] (25,1) -- (21,1) -- (21,2) -- (25,2);
    \fill[fill=black] (22,-1.5) circle (8pt);
    \draw[semithick, rounded corners=0.8pt] 
    (25,-2.5) -- (23,-2.5) -- (23,-1.5) -- (22,-1.5) -- (22,-0.5) -- (25,-0.5);
  \end{tikzpicture} 
  $$
  An almost identical
  argument entails the more general conclusion: the smallest letter in
  $w_i$ is to the left of any letter in $w_{i-1}$. Thus we can recover
  $m_i$ as the number of elements of $w_i$ that are bigger than the
  first letter of $w_i$. Consequently, the map $\psi$ is injective.
  \medskip

  {\it Part e.}  Let $\pi$ be an indecomposable avoider whose first
  letter is not $n$. To reverse $\psi$ we do as described in Part~d.
  Let $w_i$ be defined as above.  By Lemma~\ref{n-precedes-1}, $w_0$
  contains the letter $1$ and, in particular, it is nonempty. Thus, in
  the preimage of $\pi$ under $\psi$, the topmost $w_i$ will be
  nonempty, and hence $n-1$ will be to the right of $n$. By similar
  reasoning as in Part~d it also follows that the preimage avoids
  \pc\ and \pg.
\end{proof}

Given an avoider $\pi$ on $[n-1]$ and a positive integer $i$ that is
no greater than the number of left-to-right maxima in $\pi$, we define
\begin{equation} \label{phi-def}
  \begin{aligned}
    \phi(1,\pi) &= \hat\pi, \\
    \phi(i,\pi) &= \psi(\hat\pi)\;\text{ if } i>1,
  \end{aligned}
\end{equation}
where $\hat\pi$ is obtained from $\pi$ by inserting $n$ immediately
to the left of the $i$th left-to-right maximum in $\pi$.

Let $\A_n$ be the set of avoiders on $[n]$, and let $\A_n^k$ be the
subset of $\A_n$ consisting of those avoiders that have $k$
left-to-right maxima. Similarly, let $\barA_n$ be the set of
indecomposable avoiders on $[n]$, and let $\barA_n^k$ be the subset of
$\barA_n$ consisting of those indecomposable avoiders that have $k$
left-to-right maxima.

\begin{lemma}\label{lemma_phi0}
  Let $n\geq 1$. The function $\phi$, as defined by \eqref{phi-def},
  is a bijection between the Cartesian product $[k]\times\A_{n-1}^k$
  and the disjoint union $\cup_{i=1}^k\barA_{n}^i$.
\end{lemma}

\begin{proof}
  We start by showing that $\phi(i,\pi)$ has exactly $i$ left-to-right
  maxima. By definition, $\phi(i,\pi) = \psi(\hat\pi)$ where $\hat\pi$
  is obtained from $\pi$ by inserting $n$ immediately to the left of
  the $i$th left-to-right maximum in $\pi$. For the case $i=1$ we have
  $\phi(1,\pi) = \hat\pi = n\pi$ and that permutation clearly has one
  left-to-right maximum. For the case $i>1$ we note that by
  construction $\lmax\hat\pi = i$, and according to
  Lemma~\ref{lemma_psi} the number of left-to-right maxima is
  preserved under $\psi$.

  We now show that $\phi$ has the claimed codomain. The case $i=1$ is
  simple: $\phi(1,\pi)=n\pi$ has $1$ left-to-right maximum and is
  indecomposable by Lemma~\ref{n-precedes-1}. Also, $n\pi$ is an
  avoider precisely when $\pi$ is. Thus $\phi(1,\pi)$ is a member of
  $\barA_n^1$. The case $i>1$ follows from Lemma~\ref{lemma_psi0}. 

  Finally, to show that $\phi$ is bijective we give its inverse: if
  $\pi=n\tau$ then $\phi^{-1}(\pi) = (1,\tau)$; otherwise,
  $\phi^{-1}(\pi) = (\lmax\pi,\tau)$ where $\tau$ is obtained from
  $\psi^{-1}(\pi)$ by removing the largest element. That this really
  is the inverse of $\phi$ is an immediate consequence of the definition
  of $\phi$ and Lemma~\ref{lemma_psi0}.
\end{proof}

As a corollary to the preceding lemma we get that each avoider $\pi$
is of exactly one the following three forms:
\begin{itemize}
\item[] $\pi=\emptyword$,      \hfill (the empty permutation)
\item[] $\pi=\sigma\oplus \tau$,   \hfill (decomposable)
\item[] $\pi=\phi(i,\sigma)$, where $1\leq i\leq \lmax(\sigma)$,
  \hfill(indecomposable)
\end{itemize}
Note the striking similarity with the decomposition of \btrees\ as
given at the end of Section~\ref{btrees-structure}.  With avoiders as
with \btrees\ we can keep decomposing until only ``atoms'' remain.
For trees ``atom'' means the one node tree. For avoiders ``atom''
means the empty permutation. We thus get an unambiguous encoding of
avoiders; for example, the encoding of $523147896$ is
$$ 
 \phi\big(1,\phi(2,\phi(1,\emptyword)\oplus \phi(1,\emptyword))\oplus \phi(1,\emptyword)\big)\oplus
    \phi\big(3,\phi(1,\emptyword)\oplus\phi(1,\emptyword)\oplus\phi(1,\emptyword)\big).
$$
See also the example in the next section.

\section{The bijection between trees and avoiders}

Using recursion it is now easy to define a bijection---let us call
it $f$---between {\btrees} and avoiders:
\begin{equation}\label{def_f}
f(\leaf) = \emptyword,\;\;
f(\lambda_i t) = \phi_if(t),\, \text{ and }\;
f(u\oplus v) = f(u)\oplus f(v),
\end{equation}
where we have written $\lambda_i(t)$ instead of $\lambda(i,t)$, and
$\phi_i(\pi)$ instead of $\phi(i,\pi)$. Hence, to find the image of a
{\btree} $t$ under the bijection $f$, write out $t$ as a string using
$\lambda$ and $\oplus$, replace $\leaf$ with $\emptyword$ (or $\edge$
with $1$), replace $\lambda$ with $\phi$, and replace $\oplus$ (on
trees) with $\oplus$ on avoiders. Finally, translate the derived
expression to a permutation. Clearly this is an invertible process,
and so describes a bijection. For instance,
$$
\begin{tikzpicture}[ scale=0.31, baseline=-14ex ]
  \discstyle
  \exampletree
\end{tikzpicture}\,
\begin{aligned}
=\;&
  \lambda_1\Big(\lambda_2\big(\,\edge\oplus\edge\,\big)\oplus\edge\,\Big)\oplus
    \lambda_3\big(\,\edge\oplus\edge\oplus\edge\,\big)\\
\raisebox{.6ex}{$\substack{\textstyle{f}\\ \textstyle{\to}}$} \,&\,
  \phi_1\big(\phi_2(1\oplus 1)\oplus 1\big)\oplus
    \phi_3\big(1\oplus1\oplus1\big)\\
=\;& \phi_1(\phi_2(12) \oplus 1)\oplus\phi_3(123) \\
=\;& \phi_1(231\oplus 1)\oplus 2341\\
=\;& \phi_1(2314)\oplus 2341 \\
=\;& 52314\oplus 2341 = 523147896.
\end{aligned}
$$ 

\begin{remark*}
  In the course of discovering the map $\psi$, defined in
  \eqref{def_psi}, we first discovered a different map that we call
  $\theta$. The map $\theta$ can be used instead of $\psi$ in the
  above proofs with the exception that $\rpath$ on \btrees\ would not
  be mapped to $\lmin$ on avoiders; thus our main result would be
  weakened. On the other hand, $\theta$ is easier than $\psi$ to
  define: Let $\pi$ be an avoider on $[n]$ such that $n$ is not the
  first letter of $\pi$ and $n$ precedes $n-1$ (as in
  Lemma~\ref{lemma_psi0}). Also, let us write $\pi = \sigma\tau n\rho$
  in which $\tau n\rho$ is the rightmost component of $\pi$. Then
  $\theta$ is defined as in this picture:
  $$
  \begin{tikzpicture}[semithick, scale=0.5 ]
    \draw (0,0) rectangle (2,2);
    \node at (1,1) {$\sigma$};
    \draw [fill=white] (2,2) rectangle (3,3);
    \node at (2.5,2.5) {$\tau$};
    \node at (3.25,3.25) {$n$};
    \draw [fill=white] (3.5,2) rectangle (4.5,3);
    \node at (4,2.5) {$\rho$};
    
    \draw[->] (6.3,1) -- (7.8,1) node[midway, yshift=3mm] {$\theta$};

    \draw (10,0) rectangle (11,1);
    \node at (10.5,0.5) {$\tilde\tau$};
    \draw (11,1) rectangle (13,3);
    \node at (12,2) {$\tilde \sigma$};
    \node at (13.25,3.25) {$n$};
    \draw (13.5,0) rectangle (14.5,1);
    \node at (14,0.5) {$\tilde\rho$};
  \end{tikzpicture}
  $$ Here $\tilde\tau\tilde\rho$ is the standardization of $\tau\rho$
  and $\tilde\sigma$ is obtained from $\sigma$ by adding $|\tau\rho|$
  to all of the letters of $\sigma$. Note that by
  Lemma~\ref{n-precedes-1}, since $\tau n\rho$ is irreducible, the
  smallest letter of $\tau n\rho$ is found in $\rho$. Thus the letter
  $1$ in $\theta(\pi)$ belongs to $\tilde\rho$ and again using
  Lemma~\ref{n-precedes-1} we conclude that $\theta(\pi)$ is
  irreducible.
\end{remark*}

\section{Statistics on {\btrees}}\label{tree_stats}

The proofs in this section will use induction on the number of edges
in a tree.  For that reason we note that, by definition, we have
\begin{equation}\label{lemma_base_tree}
  \begin{aligned}
    &\leaves(\leaf)=1\;\text{ and}\\
    &\sub(\leaf) = \root(\leaf) = \lpath(\leaf) =
    \rpath(\leaf) = \lsub(\leaf) = \stem(\leaf)=0, 
  \end{aligned} 
\end{equation}
where $\leaf$ is the unique {\btree} with a single node. This will
serve as the basis for induction.

\begin{lemma}\label{lemma_lambda}
  The map $\lambda: [k]\times\B_{n-1}^k\to\cup_{i=1}^k\barB_{n}^i$ has
  the following properties:
  \begin{align*}
    \leaves\lambda(i,t) &= \leaves t;     \\
    \root\lambda(i,t)   &= i;             \\
    \lpath\lambda(i,t)  &= \lpath t + 1;  \\
    \rpath\lambda(i,t)  &= \rpath t + 1;  \\
    \lsub\lambda(i,t)   &= \lsub t + 1 &&\text{if } i=1; \\
    \lsub\lambda(i,t)   &= \lsub t     &&\text{if } i>1; \\
    \stemhm\lambda(i,t) &= i          &&\text{if } i\leq \stemhm t;\\
    \stemhm\lambda(i,t) &= \stemhm t  &&\text{if } i> \stemhm t.
  \end{align*} 
\end{lemma}

\begin{proof}
  Straightforward and omitted. 
\end{proof}

\begin{lemma}\label{lemma_Lambda}
  If $t=u\oplus v$, with $u\neq\leaf$ and $v\neq\leaf$, then
  \begin{align*}
    \leaves t    &= \leaves u + \leaves v;  \\
    \root t      &= \root u + \root v;      \\
    \lpath t     &= \lpath u;               \\
    \rpath t     &= \rpath v;               \\
    \lsub t      &= \lsub u;                
    \shortintertext{and, if $k$ is the largest integer such that 
      $t=\big({\oplus^k\,\edge}\,\big)\oplus v$ for some $v\neq\leaf$, then}
    \stemhm t   &= k+\stemhm v. 
  \end{align*} 
\end{lemma}

\begin{proof}
  Straightforward and omitted. 
\end{proof}

\section{Statistics on avoiders}

The proofs in this section will use induction on the number of letters
in a permutation, and we therefore note that, by definition, we have
\begin{equation}\label{lemma_base_perm}
    \asc(\emptyword) =
    \comp(\emptyword) = \lmax(\emptyword) 
    = \lmin(\emptyword) = \rmax(\emptyword) 
    = \ldr(\emptyword) = \lir(\emptyword)=0,
\end{equation}
where $\emptyword$ is the empty permutation.

\begin{lemma}\label{lemma_psi}
  The function $\psi$, as defined by \eqref{def_psi}, preserves
  left-to-right maxima, right-to-left maxima, ascents, leftmost
  decreasing run, and leftmost increasing run; in addition, it
  increases the number of left-to-right minima by one.
\end{lemma}

\begin{proof}
  Let $\pi=\sigma n \tau$, with $\sigma$ nonempty and $n-1$
  in $\tau$, be an avoider. Since
  \begin{gather*}
    \lmax\pi = 1+\lmax\sigma, \quad\;
    \rmax\pi = 1+\rmax\tau,\quad\;
    \asc\pi  = 1+\asc\sigma+\asc\tau, \\
    \ldr\pi  = \ldr \sigma,\quad\;
    \lir\pi  = \lir(\sigma n),
  \end{gather*}
  and $\psi$ preserves the patterns of $\sigma$ and $\tau$, it
  immediately follows that $\psi$ preserves $\lmax$, $\rmax$, $\asc$,
  $\ldr$, and $\lir$.  

  The reason why $\psi$ increases $\lmin$ by one is a bit more involved:
  Note first that, with the same notation as in the definition of
  $\psi$, we have
  $$\lmin\pi = \lmin\sigma_1 + \lmin w_0.
  $$ Moreover, any left-to-right minimum in $w_0$ is to the right of
  the smallest element of $w_1$. This is because otherwise we have a
  $\pg$ pattern $bdac$ in which $c=\min w_1$; $a$ is a left-to-right
  minimum to the left of $c$; $d$ is the element immediately to the
  left of $a$ in $\pi$; and $b$ is any element in $\sigma_1$. (Here,
  the assumption that $n-1$ is in $\tau$ entails that $w_1$ is
  nonempty and the assumption that $\tau$ is nonempty entails that
  $\sigma_1$ is nonempty.) The subsequence $w_0$ of $\tau$ is fixed
  under $\psi$. Also, $\psi$ preserves the pattern of $\sigma_1$
  but adds $m_i+1$ to each of its elements. Exactly one new
  left-to-right minimum is thus obtained, namely $1+\max w_0$, the
  image of $\min w_1$ under $\psi$.
\end{proof}

\begin{lemma}\label{lemma_phi}
  The function $\phi$, as defined by \eqref{phi-def}, has the
  following properties.
  \begin{align}
    \asc\phi(i,\pi)  &= \asc\pi       \label{asc.phi};\\
    \lmax\phi(i,\pi) &= i             \label{lmax.phi}; \\
    \lmin\phi(i,\pi) &= \lmin\pi + 1  \label{lmin.phi};\\
    \rmax\phi(i,\pi) &= \rmax\pi + 1  \label{rmax.phi}; \\
    \ldr\phi(i,\pi)  &= \ldr\pi + 1 &&\text{if } i=1; \label{ldr.phi_a}\\
    \ldr\phi(i,\pi)  &= \ldr\pi     &&\text{if } i>1; \label{ldr.phi_b}\\
    \lir\phi(i,\pi)  &= i        &&\text{if } i\leq\lir\pi; \label{lir.phi_a}\\
    \lir\phi(i,\pi)  &= \lir\pi  &&\text{if } i>\lir\pi.    \label{lir.phi_b}
  \end{align}
\end{lemma}

\begin{proof}
  For the case $i=1$ we have $\phi(1,\pi) = \hat\pi = n\pi$, and all
  of the statements \eqref{asc.phi} through \eqref{lir.phi_b} follow
  immediately. The interesting case is $i>1$, which we now consider.
  By definition, $\phi(i,\pi) = \psi(\hat\pi)$ where $\hat\pi$ is
  obtained from $\pi$ by inserting $n$ immediately to the left of the
  $i$th left-to-right maximum in $\pi$. By construction $\lmax\hat\pi
  = i$, and according to Lemma~\ref{lemma_psi} the number of
  left-to-right maxima is preserved under $\psi$; thus
  $\lmax\phi(i,\pi) = i$, proving \eqref{lmax.phi}. Further,
  $$ 
  \asc\hat\pi = \asc\pi; \quad
  \lmin\hat\pi = \lmin\pi; \quad
  \rmax\hat\pi = \rmax\pi + 1; \quad
  \ldr\hat\pi = \ldr\pi.
  $$ Hence statements \eqref{asc.phi}, \eqref{lmin.phi},
  \eqref{rmax.phi}, and \eqref{ldr.phi_b} follow from the
  corresponding statements in Lemma~\ref{lemma_psi}. For
  \eqref{lir.phi_a} and \eqref{lir.phi_b} note that an element in the
  leftmost increasing run is a left-to-right maximum, and that
  inserting $n$ in front of the $i$th of those elements (in order to
  build $\hat\pi$) results in $\ldr\hat\pi = i$. On the other hand,
  inserting $n$ in front of a left-to-right maximum that is not in the
  leftmost increasing run results in $\ldr\hat\pi = \lir\pi$. Thus,
  also, \eqref{lir.phi_a} and \eqref{lir.phi_b} follow from the
  corresponding statements in Lemma~\ref{lemma_psi}.
\end{proof}

\begin{lemma}\label{lemma_prod}
  With $\pi=\sigma\oplus\tau$ we have
  \begin{align*}
    \comp\pi &= \comp\sigma + \comp\tau;  && \\
    \asc\pi  &= 1+\asc\sigma+\asc\tau;    && \\
    \lmax\pi &= \lmax\sigma + \lmax\tau;  && \\
    \lmin\pi &= \lmin\sigma;              && \\
    \rmax\pi &= \rmax\tau;                && \\
    \ldr\pi  &= \ldr\sigma;               && \\
    \shortintertext{and, if $k$ is the largest integer such that 
      $\pi=(\oplus^k 1)\oplus \tau$ for some nonempty $\tau$, then}
    \lir\pi  &= k + \lir\tau              && 
  \end{align*}
\end{lemma}

\begin{proof}
  Straightforward and omitted. 
\end{proof}

\begin{theorem}\label{thm_f}
  Let $f$ be the bijection from {\btrees} on $n+1$ nodes onto
  length $n$ avoiders, as defined by \eqref{def_f}. It sends the first
  7-tuple of statistics, below, to the second 7-tuple.
  $$
  \begin{array}{llllllllll}
    ( &\sub,  &\leaves, &\root, &\lpath, &\rpath, &\lsub, &\stemhm & )\\
    ( &\comp, &1+\asc,  &\lmax, &\lmin,  &\rmax,  &\ldr,  &\lir    & )
  \end{array}
  $$
\end{theorem}

\begin{proof}
  The proof proceeds by induction. The base
  case follows from Lemma~\ref{lemma_base_perm} and
  Lemma~\ref{lemma_base_tree}. Let $t$ be any {\btree} with $n+1$
  nodes. We split into two cases: (1) $t=\lambda_i u$ is
  indecomposable; (2) $t=u\oplus v$ is decomposable. 
  
  Case 1: That $\comp f(\lambda_is) = \comp\phi_if(s) = 1 =
  \sub\lambda_is$ is clear since both $\phi_i\pi$ and $\lambda_it$ are
  indecomposable for any permutation $\pi$ and any {\btree} $t$. Also,
  \begin{align*}
    1+\asc f(\lambda_is) 
    &= 1+\asc\phi_if(s)   &&\text{by definition of $f$}\\
    &= 1+\asc f(s)        &&\text{by Lemma~\ref{lemma_phi}}\\
    &= \leaves s          &&\text{by induction}\\
    &= \leaves\lambda_is  &&\text{by Lemma~\ref{lemma_lambda}}
  \end{align*}
  Thus we have proved that, for indecomposable trees, $f$ sends $\sub$
  to $\comp$ and $\leaves$ to $1+\asc$. The proofs of the remaining
  statements follow the same pattern: recall the definition of $f$,
  apply Lemma~\ref{lemma_phi}, use the induction hypothesis, and
  finish by applying Lemma~\ref{lemma_lambda}.

  Case 2: We have 
  \begin{align*}
    \comp f(u\oplus v) 
    &= \comp (f(u)\oplus f(v))  &&\text{by definition of $f$}\\
    &= \comp f(u) + \comp f(v)  &&\text{by Lemma~\ref{lemma_prod}}\\
    &= \sub u + \sub v          &&\text{by induction}\\
    &= \sub(u\oplus v)          &&\text{by Lemma~\ref{lemma_Lambda}}
  \end{align*} 
  Again, the remaining statements follow similarly.
\end{proof}

\section{An involution on {\btrees}}\label{h}

In this section we define an involution on {\btrees}. To that end we
now describe a new way of decomposing {\btrees}. Schematically the
sum $\oplus$ on {\btrees} is described by
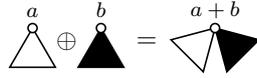
\begin{center}
  \makebox[3.3cm][l]{
  \begin{tikzpicture}[semithick, scale=0.3]
    \def\tri{ -- +(1,-1.732) -- +(-1,-1.732) -- cycle};
    \draw (0,0) \tri;
    \filldraw[fill=white] (0,0) circle (2mm);
    \node[font=\footnotesize] at (0,0.7) {$a$};
    \node at (1.5,-0.6) {$\oplus$};
    \filldraw (3,0) \tri;
    \filldraw[fill=white] (3,0) circle (2mm);
    \node[font=\footnotesize] at (3,0.8) {$b$};
    \node at (5,-0.6) {=};
    \draw[shift={(8,0)}] (0,0) [rotate=-45] \tri;
    \filldraw[shift={(8,0)}] (0,0) [rotate=45 ] \tri;
    \filldraw[fill=white] (8,0) circle (2mm);
    \node[font=\footnotesize] at (8,0.8) {$a+b$};
  \end{tikzpicture}
  }
\end{center}
An alternative sum is\\[-5ex]
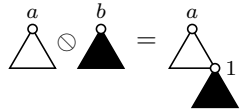
\begin{center}
  \makebox[3.3cm][l]{
  \begin{tikzpicture}[semithick, scale=0.3]
    \def\tri{ -- +(1,-1.732) -- +(-1,-1.732) -- cycle};
    \draw (0,0) \tri;
    \filldraw[fill=white] (0,0) circle (2mm);
    \node[font=\footnotesize] at (0,0.7) {$a$};
    \node at (1.5,-0.6) {$\obslash$};
    \filldraw (3,0) \tri;
    \filldraw[fill=white] (3,0) circle (2mm);
    \node[font=\footnotesize] at (3,0.8) {$b$};
    \node at (5,-0.6) {=};
    \draw (7,0) \tri;
    \filldraw (8,-1.732) \tri;
    \filldraw[fill=white] (7,0) circle (2mm);
    \filldraw[fill=white] (8,-1.732) circle (2mm);
    \node[font=\footnotesize] at (7,0.7) {$a$};
    \node[font=\footnotesize] at (8.7,-1.732) {$1$};
  \end{tikzpicture}
  }
\end{center}
That is, to get $u\obslash v$ we join $u$ and $v$ by identifying the
rightmost leaf in $u$ with the root of $v$, and that node is assigned
the label $1$. Note that
\begin{align}
  \root(u\oplus v)  &= \root(u) + \root(v)\\
  \rpath(u\oplus v) &= \rpath(v) \\
  \shortintertext{while}
  \root(u\obslash v)  &= \root(u)\\
  \rpath(u\obslash v) &= \rpath(u) + \rpath(v) \label{rpath_obslash}.
\end{align}
for $u\neq\leaf$ and $v\neq\leaf$. Thus, with respect to $\obslash$,
$\rpath$ plays the role of $\root$, and vice versa. There is also a
map $\gamma$ that plays a role analogous to that of $\lambda$:
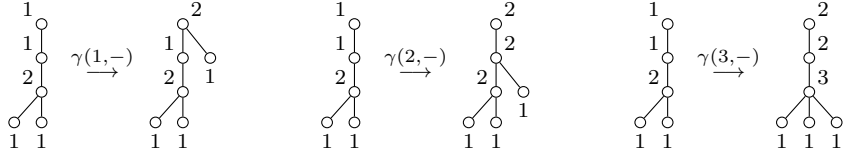
$$
\begin{tikzpicture}[scale=\scl, baseline=(r11.base) ]
  \eeev;
  \draw (r) node[above left=-0.5pt] {1} -- 
  (r1) node[above left=-0.5pt] {1} -- 
  (r11) node[above left=-0.5pt] {2};
\end{tikzpicture}
\;\;\raisebox{2ex}{$\substack{\gamma(1,\slot)\\ \lra}$}\;
\begin{tikzpicture}[ scale=\scl, baseline=(r11.base) ]
  \eeev;
  \node [disc] (r2) at ( 0.8, 2 ) {};
  \draw (r) node[above right=-0.5pt] {2} -- 
  (r1) node[above left=-0.5pt] {1} -- 
  (r11) node[above left=-0.5pt] {2};
  \draw (r) -- (r2) node[below=1pt] {1};
\end{tikzpicture}
\qquad\quad\;
\begin{tikzpicture}[scale=\scl, baseline=(r11.base) ]
  \eeev;
  \draw (r) node[above left=-0.5pt] {1} -- 
  (r1) node[above left=-0.5pt] {1} -- 
  (r11) node[above left=-0.5pt] {2};
\end{tikzpicture}
\;\;\raisebox{2ex}{$\substack{\gamma(2,\slot)\\ \lra}$}\;
\begin{tikzpicture}[ scale=\scl, baseline=(r11.base) ]
  \eeev;
  \node [disc] (r12) at ( 0.8, 1 ) {};
  \draw (r) node[above right=-0.5pt] {2} -- 
  (r1) node[above right=-0.5pt] {2} -- 
  (r11) node[above left=-0.5pt] {2};
  \draw (r1) -- (r12) node[below=1pt] {1};
\end{tikzpicture}
\qquad\quad\;
\begin{tikzpicture}[scale=\scl, baseline=(r11.base) ]
  \eeev;
  \draw (r) node[above left=-0.5pt] {1} -- 
  (r1) node[above left=-0.5pt] {1} -- 
  (r11) node[above left=-0.5pt] {2};
\end{tikzpicture}
\;\;\raisebox{2ex}{$\substack{\gamma(3,\slot)\\ \lra}$}\;
\begin{tikzpicture}[ scale=\scl, baseline=(r11.base) ]
  \eeev;
  \node [disc] (r112) at ( 0.8, 0.1 ) {};
  \draw (r) node[above right=-0.5pt] {2} -- 
  (r1) node[above right=-0.5pt] {2} -- 
  (r11) node[above right=-0.5pt] {3};
  \draw (r11) -- (r112) node[below=1pt] {1};
\end{tikzpicture}
$$ Here is how $\gamma(i,t)$ is defined in general: Assume that the
length of the right path of $t$ is $k$ and that $i$ is an integer such
that $1\leq i\leq k$. Let us by $x$ refer to the $i$th node on the
right path of $t$. Then $\gamma(i,t)$ is obtained from $t$ by joining
a new leaf via an edge to $x$, making the new leaf the rightmost leaf
in $\gamma(i,t)$; and, lastly, adding $1$ to the label of each node on
the new right path, except for the new leaf (in which the new right
path ends). Note, in particular, that $\rpath\gamma(i,t) = i$.

We now connect the two ways we have to decompose {\btrees} by
defining an endofunction $h:\B\to\B$:
$$h(\leaf) = \leaf,\;\;
h(\lambda_it) = \gamma_ih(t),\,\text{ and }\; 
h(u\oplus v) = h(v) \obslash h(u).
$$
For instance,
$$
\begin{tikzpicture}[ scale=0.35, baseline=-7.5ex ]
  \discstyle
  \exampletree
\end{tikzpicture}\,
\begin{aligned}
\;=\;\;&
  \lambda_1\Big(\lambda_2\big(\,\edge\oplus\edge\,\big)\oplus\edge\,\Big)\oplus
    \lambda_3\big(\,\edge\oplus\edge\oplus\edge\,\big)\\[0.9ex]
\raisebox{.6ex}{\;\;$\substack{\textstyle{h}\\ \textstyle{\to}}$}\;&
\gamma_3\big(\,\edge\obslash\edge\obslash\edge\,\big)
\obslash\gamma_1\Big(\,\edge\obslash\gamma_2\big(\,\edge\obslash\edge\,\big)\Big)
\;=\,
\end{aligned}
\begin{tikzpicture}[scale=0.35, baseline=(r1.base) ]
  \style
  \node [disc] (r)       at ( 0,   7.4 ) {};
  \node [disc] (r1)      at ( 0,   6.2 ) {};
  \node [disc] (r11)     at ( 0,   5   ) {};
  \node [disc] (r111)    at (-0.8, 4   ) {};
  \node [disc] (r112)    at ( 0.8, 4   ) {};
  \node [disc] (r1121)   at ( 0,   3   ) {};
  \node [disc] (r11211)  at ( 0,   1.8 ) {};
  \node [disc] (r112111) at (-0.8, 0.8 ) {};
  \node [disc] (r112112) at ( 0.8, 0.8 ) {};
  \node [disc] (r1122)   at ( 1.6, 3   ) {};
  \draw 
  (r) node[above=1pt] {2} -- (r1) node[above right=-1pt] {2} 
  (r1) -- (r11) node[above right=-.8pt] {2} 
  (r11) -- (r111) node[left=1pt] {1}
  (r11) -- (r112) node[above right=-.8pt] {1}
  (r112) -- (r1121) node[left=1pt] {1}
  (r1121) -- (r11211) node[above right=-.8pt] {2}
  (r11211) -- (r112111) node[below=1pt] {1}
  (r11211) -- (r112112) node[below=1pt] {1}
  (r112) -- (r1122) node[below=1pt] {1};
\end{tikzpicture}
$$ 
It should be clear that $h$ is defined to translate an encoding
based on $\lambda$ and $\oplus$ into an encoding based on $\gamma$ and
$\obslash$.  As it turns out, $h$ is an involution! For the proof of
the following theorem we refer the reader to a forthcoming paper by
the present authors \cite{CKS}, but first we need to define
the statistic $\stemh(t)$.  This is simply the statistic $\stemhm$ on
the mirror image of $t$, where the mirror image $m$ is the involution
on {\btrees} that recursively reverses the order of subtrees (see the
end of Section~\ref{prel} for the definition of $\stemhm$).  To be
precise, we have $m(\leaf)=\leaf$, $m(\lambda_it)=\lambda_im(t)$, and
$m(u\oplus v)= m(v)\oplus m(u)$.  Then we have
$\stemh(t)=\stemhm(m(t))$.

Another way to define $\stemh(t)$ is as follows: Using $\gamma$ and
$\obslash$ we can write
$t=\gamma_{i_1}(\gamma_{i_2}(\dots\gamma_{i_k}(u)))$ in which $u$ is
either a single node or decomposable with respect to $\obslash$, and
we then let $\stemh(t)=k$.

\begin{theorem}\label{thm_h}
  On {\btrees} with at least one edge, the function $h$ is an
  involution, and it sends the first tuple below to the second.
  $$
  \begin{array}{lllllllllll}
    ( &\leaves,    &\internal, &\root,  &\rpath, &\sub,  &\rsub, &\stem,  &\stemh & ) \\
    ( & \internal, &\leaves,   &\rpath, &\root,  &\rsub, &\sub,  &\stemh, &\stem  & )
  \end{array}
  $$
\end{theorem}

\begin{corollary}\label{cor_hf}
  On length $n$ avoiders, the involution $f^{-1}\circ h\circ f$ sends 
  $$
  (\,\asc,\,\lmax,\,\rmax\,)\;\text{ to }\; 
  (\,\des,\,\rmax,\,\lmax\,).
  $$
\end{corollary}
\begin{proof}
  Follows from combining Theorem~\ref{thm_h} with Theorem~\ref{thm_f}.
\end{proof}

\begin{corollary}\label{cor_hmf}
  On length $n$ avoiders, the involution $f^{-1}\circ m \circ h\circ m
  \circ f$ sends
  $$
  (\,\asc,\,\lmax,\,\lmin,\,\comp,\,\ldr\,)\;\text{ to }\; 
  (\,\des,\,\lmin,\,\lmax,\,\ldr,\,\comp\,).
  $$
\end{corollary}
\begin{proof}
  Follows from Theorems \ref{thm_h} and \ref{thm_f} together with the
  definition of $m$.
\end{proof}

We end this section with the observation that $h$ restricted to
\btrees\ with all nodes labeled 1 (except the root) induces an
involution on unlabeled rooted plane trees. This involution appears to
be new and its consequences will be explored in a forthcoming paper
\cite{CKS}. In particular, this yields results akin to
Corollaries~\ref{cor_hf} and \ref{cor_hmf} for one-stack sortable
permutations.  Moreover, this also gives rise to a genuinely new
bijection between (1-2-3)-avoiding and (1-3-2)-avoiding permutations
and yields new equidistributions of statistics on these two classes of
permutations. (See \cite{CK} for a classification of the known
bijections between (1-2-3)-avoiding and (1-3-2)-avoiding
permutations.)

\section{A conjecture about two-stack sortable permutations}

Dulucq et al~\cite{DuGiGu98} proved that the pair $(\asc,\lmax)$ on
avoiders is equidistributed with the pair $(\des,\rmax)$ on $2$-stack
sortable permutations. We make the following conjecture:

\begin{conjecture}
  The quadruple
  $ (\,\comp,\,\asc,\,\ldr,\,\rmax\,)$ 
  has the same distribution on length $n$ avoiders as it has on $2$-stack
  sortable permutations of length $n$.
\end{conjecture}

\section{Acknowledgments}

The first author did some of the work on this paper while visiting the
Mittag-Leffler Institute, and he is grateful for their hospitality.

\bibliographystyle{plain}

\end{document}